\documentclass[11pt,a4paper]{article}
\usepackage[T1]{fontenc}
\usepackage{lmodern}
\usepackage[margin=1in]{geometry}
\usepackage{amsmath,amssymb,amsthm}
\usepackage{microtype}
\usepackage{enumitem}
\usepackage{cite}
\usepackage[hidelinks]{hyperref}
\hypersetup{pdftitle={Understanding PAGE: Note III on Recursive Variance Reduction}}
\newtheorem{assumption}{Assumption}
\newtheorem{lemma}{Lemma}
\newtheorem{theorem}{Theorem}
\newcommand{\E}{\mathbb{E}}
\newcommand{\R}{\mathbb{R}}
\newcommand{\norm}[1]{\left\|#1\right\|}
\newcommand{\ip}[2]{\left\langle #1,#2\right\rangle}
\newcommand{\eps}{\varepsilon}
\setlist[itemize]{leftmargin=1.4em,itemsep=2pt,topsep=4pt,parsep=0pt}
\title{Variance Reduction Note 3: The Convergence of PAGE}
\author{Alina Lodhi\\
xieyunqi27@gmail.com}
\date{}

\begin{document}
\maketitle
\vspace{-2em}
\begin{abstract}
The probabilistic gradient estimator (PAGE) alternates randomly between
an accurate gradient refresh and a cheap recursive correction. This third
expository note on recursive variance reduction explains the algorithm
and gives a self-contained convergence proof for smooth nonconvex
optimization. A Lyapunov function combines objective decrease with
estimation error, turning the refresh probability into an explicit
contraction parameter. The resulting expected oracle complexity is
$O(n+\sqrt n\,\eps^{-2})$ for finite sums and $O(\eps^{-3})$ for
online problems, with the other problem parameters fixed. Initialization,
refreshes, and both evaluations in each gradient difference are counted.
The presentation is pedagogical and makes no claim of a new algorithm
or convergence rate.
\end{abstract}

\section{Random refreshes and recursive corrections}
Consider the smooth, possibly nonconvex objective
\begin{equation}\label{eq:problem}
 \min_{x\in\R^d}F(x):=\E_\xi[f(x;\xi)],\qquad
 g(x;\xi):=\nabla f(x;\xi).
\end{equation}
Recursive estimators such as SARAH and SPIDER
\cite{sarah2017,spider2018} use shared-sample differences to track
changes in the gradient. A refresh prevents the accumulated tracking
error from becoming too large. PAGE, introduced by Li et al.\
\cite{page2021}, performs this refresh with a specified probability
at every iteration. The resulting random refresh times remove the need
for a prescribed epoch length.

The method has two computational scales: an expensive refresh batch and
a smaller batch of gradient differences. The refresh probability balances
their expected costs. Its mathematical role is equally direct: refreshing
replaces the old estimation error, so that error survives with probability
$1-p$. This yields a short Lyapunov proof with a constant step size.
Related probabilistic designs include loopless SVRG
\cite{loopless2019}; STORM instead incorporates a damped recursive
correction at every step \cite{storm2019}. For a general introduction to
variance reduction, see \cite{gower2020}.

\section{The oracle and the PAGE update}
\begin{assumption}\label{ass:oracle}
The differentiable function $F$ is bounded below by a finite number
$F_\star$. For some $L>0$ and all $x,y\in\R^d$,
\begin{equation}\label{eq:oracle}
 \E[g(x;\xi)]=\nabla F(x),\qquad
 \E\norm{g(x;\xi)-g(y;\xi)}^2\leq L^2\norm{x-y}^2.
\end{equation}
For online problems, additionally assume, for some $\sigma\geq0$,
\begin{equation}\label{eq:variance}
 \E\norm{g(x;\xi)-\nabla F(x)}^2\leq\sigma^2.
\end{equation}
Fresh samples are independent of the past and of one another. A sampled
function can be evaluated at two specified points. Norms are Euclidean.
\end{assumption}
Jensen's inequality implies that $F$ is $L$-smooth. For a finite sum
$F(x)=n^{-1}\sum_{i=1}^n f_i(x)$, use uniform indices and
$g(x;i)=\nabla f_i(x)$. Exact full-gradient refreshes then eliminate
the need for \eqref{eq:variance}.

Choose a fixed step size $\eta>0$, a probability $p\in(0,1]$, and
integer batch sizes $B,b\geq1$. Initialize $v_0$ with a refresh at
a deterministic $x_0$. For $t=0,\ldots,T-1$, compute
\begin{equation}\label{eq:step}
 x_{t+1}=x_t-\eta v_t,
\end{equation}
then draw a fresh independent Bernoulli variable $C_{t+1}$ with
$\Pr(C_{t+1}=1)=p$, and set
\begin{equation}\label{eq:page}
 v_{t+1}=\begin{cases}
 \displaystyle\frac1B\sum_{i=1}^B g(x_{t+1};\zeta_{t+1,i}),
       &C_{t+1}=1,\\[1.2ex]
 \displaystyle v_t+\frac1b\sum_{i=1}^b
 [g(x_{t+1};\xi_{t+1,i})-g(x_t;\xi_{t+1,i})],
       &C_{t+1}=0.
 \end{cases}
\end{equation}
Both evaluations in each difference use the same sample. For finite sums,
replace every refresh, including initialization, by the full gradient;
set $B=n$ for cost accounting and sample inner indices independently
with replacement. Return $x_R$ for an independent uniform
$R\in\{0,\ldots,T-1\}$. The expected cost of this implementation is
\begin{equation}\label{eq:cost}
 \E N=B+T[pB+2(1-p)b].
\end{equation}
This count includes the final estimator update, which can be omitted
without changing the returned point or the guarantee.

\section{A short convergence proof}
Write $e_t=v_t-\nabla F(x_t)$, $Q_t=\E\norm{e_t}^2$,
$D_t=\E\norm{v_t}^2$, and $G_t=\E\norm{\nabla F(x_t)}^2$.
Let $\nu^2=\sigma^2/B$ online and $\nu^2=0$ for exact finite-sum
refreshes. In both settings $Q_0\leq\nu^2$.

\begin{lemma}[Tracking error]\label{lem:error}
Under Assumption~\ref{ass:oracle},
\begin{equation}\label{eq:error}
 Q_{t+1}\leq(1-p)Q_t+
 \frac{(1-p)L^2\eta^2}{b}D_t+p\nu^2.
\end{equation}
\end{lemma}
\begin{proof}
Condition on all randomness used to construct $x_t,v_t$; then
$x_{t+1}$ is also fixed. A refresh has conditional mean-square error
at most $\nu^2$. In the recursive branch, write
$e_{t+1}=e_t+Z_{t+1}$, where
\[
 Z_{t+1}=\frac1b\sum_{i=1}^b
 [g(x_{t+1};\xi_i)-g(x_t;\xi_i)]
 -[\nabla F(x_{t+1})-\nabla F(x_t)].
\]
The increment has conditional mean zero. Independence within the batch,
the variance identity, and \eqref{eq:oracle} bound its conditional
second moment by $L^2\norm{x_{t+1}-x_t}^2/b$.
Its cross term with $e_t$ vanishes. Average over the fresh Bernoulli
choice and the past, and use \eqref{eq:step}.
\end{proof}

\begin{lemma}[Descent]\label{lem:descent}
For the update \eqref{eq:step},
\begin{equation}\label{eq:descent}
 \E F(x_{t+1})\leq\E F(x_t)
 -\frac\eta2G_t-\frac\eta2(1-L\eta)D_t+\frac\eta2Q_t.
\end{equation}
\end{lemma}
\begin{proof}
Apply $L$-smoothness to $x_{t+1}=x_t-\eta v_t$ and substitute
\[
 2\ip{\nabla F(x_t)}{v_t}
 =\norm{\nabla F(x_t)}^2+\norm{v_t}^2-\norm{e_t}^2.
\]
Taking expectations gives \eqref{eq:descent}; conditional unbiasedness
of $v_t$ is not required.
\end{proof}

\begin{theorem}[Expected stationarity]\label{thm:main}
Let $\Delta=F(x_0)-F_\star$ and $c=(1-p)/(pb)$. If
\begin{equation}\label{eq:stepsize}
 0<\eta\leq\frac1{L(1+\sqrt c)},
\end{equation}
then PAGE satisfies
\begin{equation}\label{eq:rate}
 \E\norm{\nabla F(x_R)}^2
 \leq\frac{2\Delta}{\eta T}+\frac{Q_0}{pT}+\nu^2
 \leq\frac{2\Delta}{\eta T}+\frac{\nu^2}{pT}+\nu^2.
\end{equation}
\end{theorem}
\begin{proof}
Use the nonnegative Lyapunov quantity
$\Phi_t=\E[F(x_t)-F_\star]+\eta Q_t/(2p)$.
Multiply \eqref{eq:error} by $\eta/(2p)$ and add
\eqref{eq:descent}. The $Q_t$ terms cancel, giving
\[
 \Phi_{t+1}-\Phi_t\leq-\frac\eta2G_t
 -\frac\eta2[1-L\eta-cL^2\eta^2]D_t
 +\frac\eta2\nu^2.
\]
The bracket is nonnegative: with $r=\sqrt c$, condition
\eqref{eq:stepsize} gives
\[
 L\eta+cL^2\eta^2
 \leq\frac1{1+r}+\frac{r^2}{(1+r)^2}\leq1.
\]
Sum over $t=0,\ldots,T-1$ and use $\Phi_T\geq0$ and
$\Phi_0=\Delta+\eta Q_0/(2p)$. Division by $\eta T/2$ yields
\eqref{eq:rate}, because $R$ is independent and uniform.
\end{proof}

\paragraph{Finite-sum complexity.}
Take exact refreshes, $B=n$, $b=\lceil\sqrt n\rceil$,
$p=b/(n+b)$, and equality in \eqref{eq:stepsize}. Then
$c=n/b^2\leq1$, so $\eta\geq1/(2L)$, and $Q_0=\nu^2=0$.
For a target $\eps>0$, choosing
$T=\max\{1,\lceil4L\Delta/\eps^2\rceil\}$ gives
$\E\norm{\nabla F(x_R)}^2\leq\eps^2$.
Moreover, \eqref{eq:cost} gives
\begin{equation}\label{eq:finite-cost}
 \E N=n+\frac{3nb}{n+b}T
 \leq n+3bT
 =O\!\left(n+\frac{\sqrt n\,L\Delta}{\eps^2}\right).
\end{equation}
The $O(\sqrt n)$ cost from rounding $T$ is absorbed by $n$.

\paragraph{Online complexity.}
Choose
\begin{equation}\label{eq:online}
 B=\max\{1,\lceil2\sigma^2/\eps^2\rceil\},\qquad
 b=\lceil\sqrt B\rceil,\qquad p=\frac b{B+b},
\end{equation}
and equality in \eqref{eq:stepsize}. Again $c=B/b^2\leq1$,
while $\nu^2\leq\eps^2/2$. With
\begin{equation}\label{eq:horizon}
 T=\left\lceil\frac{4\Delta}{\eta\eps^2}+\frac1p\right\rceil,
\end{equation}
the first two terms in \eqref{eq:rate} sum to at most $\eps^2/2$.
Thus $\E\norm{\nabla F(x_R)}^2\leq\eps^2$. Since
$b/p=B+b$ and $\eta^{-1}\leq2L$,
\begin{align}\label{eq:online-cost}
 \E N&\leq B+3bT
 =O\!\left(B+\frac{\sqrt B\,L\Delta}{\eps^2}\right)\notag\\
 &=O\!\left(1+\frac{\sigma^2}{\eps^2}
       +\left(1+\frac\sigma\eps\right)\frac{L\Delta}{\eps^2}\right).
\end{align}
For fixed problem parameters this is $O(\eps^{-3})$ as
$\eps\downarrow0$. In either setting, Cauchy--Schwarz also gives
$\E\norm{\nabla F(x_R)}\leq\eps$.

\paragraph{Interpretation and scope.}
The term $Q_0/(pT)$ shows why initialization belongs in the proof.
A sampled refresh leaves a noise floor $\nu^2$, which must decrease
with the requested accuracy. The probability $p=b/(B+b)$ and
$b\asymp\sqrt B$ balance the expected work while permitting a step
size of order $1/L$. Cost bounds are in expectation over the random
refreshes. The chosen horizons use $\Delta$ or an upper bound on it;
online batch selection also uses a variance bound. The theorem concerns
a random iterate and first-order stationarity. Linear rates under a
Polyak--\L{}ojasiewicz condition need an additional assumption and
are treated in the original PAGE paper \cite{page2021}.

\section{Recommended reading}

\paragraph{Variance-reduction methods.}
Basic guarantees \cite{ghadimi2013,gower2020,unifiedsgd2020}; gradient memories and control
variates \cite{sag2012,sdca2013,svrg2013,s2gd2013,saga2014,finito2014,proxsvrg2014}; sampling and acceleration \cite{catalyst2015,katyusha2016,svrgpp2016,katyushax2018,sagasampling2019,jacsketch2018,svrgsd2018,loopless2019,looplesssampling2021}; nonconvex and recursive
estimators \cite{reddi2016,allenzhu2016,natasha2017,sarah2017,isarah2018,scsg2017,spider2018,snvrg2018,spiderboost2019,proxsarah2020,zerosarah2021,page2021,nonconvexsampling2019}; momentum and adaptive schedules \cite{storm2019,hybrid2019,stormplus2021,superadam2021,metastorm2022,adaspider2022,jiang_adaptive2024,adasvrg2021,almosttune2020,adaptiveextra2022,alfpg2026};
second-order extensions \cite{ssrgd2019,natasha22018,svrc2019,svrcnewton2018}; lower bounds \cite{lower2023}.

\paragraph{Compositional optimization.}
Two-level  approximation \cite{scgd2014,asc2017,nasa2020,nearlyeasy2020}; finite-sum and recursive
estimators \cite{lian2017,compavg2018,civr2019,sarahcomp2019}; multilevel tracking \cite{nested2019,jiang_multilevel2022,levelindependent2020,jiang_revisiting2025};
coupled and partially observed blocks \cite{fccotheory2022,jiang_msvr2022,weakcomp2023,jiang_fcco2026,momentumcomp2025}; projection-free
methods \cite{vrfw2016,projectioncomp2022,jiang_projection2024,scgpp2019,jiang_projection2026}; statistical-model applications and federated
extensions \cite{jiang_models2023,fednest2022,jiang_models2026,localcomp2022}.

\paragraph{Sign-based optimization.}
Quantization, signSGD, majority voting, and error feedback \cite{qsgd2017,terngrad2017,signsgd2018,majority2019,errorfeedback2019};
variance-reduced sign updates and Lion \cite{jiang_sign2024,lion2023,jiang_lion2026}; compressed gradient
tracking \cite{diana2019,marina2021,jiang_distributed2026,ef212021,dasha2022}. Related distributed control-variate methods include
SCAFFOLD, FedPAGE, decentralized hybrid estimators, and ProxSkip \cite{scaffold2020,fedpage2021,hybriddecentral2021,proxskip2022}.

\paragraph{Minimax optimization.}
Finite-sum saddle-point methods, matrix games, and variance-reduced
extragradient \cite{saddle2016,matrixgames2019,svre2019,catalystminimax2020}; recursive, hybrid, adaptive, and decentralized
nonconvex minimax methods \cite{sreda2020,hybridminimax2020,adgda2023,dream2024}; a variational-inequality treatment \cite{vi2022}.


\end{document}